\documentclass[12pt]{article}
\input{epsf}
\newtheorem{theorem}{Theorem}
\newtheorem{lemma}{Lemma}

\newenvironment{proof}{{\bf Proof.}}{\hfill\rule{2mm}{2mm}}

\newtheorem{remarka}{Remark}

\def\lcs#1{{\rm lcs}(#1)}
\def\scs#1{{\rm scs}(#1)}
\def\max {{\rm \  max \ }}
\def\per {{\rm per}}

\title{\bf On the size of the minimum critical set of a Latin square}
\author{
{\bf M. Ghandehari$^a$, H. Hatami$^b$ and
E.S. Mahmoodian$^a$} \\
$^a${\small\it Department of Mathematical Sciences}\\
$^b${\small\it Department of Computer Engineering }\\
{\small Sharif University of Technology} \\
{\small P.O. Box 11365--9415, Tehran, I.R. Iran}}
\date{}
\begin{document}
\maketitle
\begin{abstract}

A {\sf critical set } in an $n \times n$ array is a set $C$ of
given entries, such that there exists a unique extension of $C$ to
an $n\times n$ Latin square and no proper subset of $C$ has this
property. For a Latin square $L$,  $\scs{L}$
 denotes the size of the
{\sf smallest critical set} of $L$, and $\scs{n}$ is the minimum
of $\scs{L}$ over all Latin squares $L$ of order $n$. We find an
upper bound for the number of partial Latin squares of size $k$
and prove that $$n^2-(e+o(1))n^{10/6} \le \max \scs{L} \le
n^2-\frac{\sqrt{\pi}}{2}n^{9/6}.$$
This improves a result of N. Cavenagh (Ph.D. thesis, The
University of Queensland, 2003) and disproves one of his
conjectures. Also it improves the previously known lower bound for
the size of the largest critical set of any Latin square of order
$n$.
\end{abstract}
{{\sc Keywords:} Critical sets; Latin squares; Partial Latin
Squares.
%%%%%%%%%%%%%%%%%%%%%%%%%%%%%%%%%%%%%%%%%%%%%%%%%%%%%%%%%%%%%%%%%%%%%

\section{Introduction}        %sec1
A {\sf Latin square} of order $n$ is an $n \times n$ array of
integers, chosen from the set $X = \{1,2, \ldots, n\}$ such that
each element of $X$ occurs exactly once in each row and exactly
once in each column. A Latin square can also be written as a set
of ordered triples $\{ (i,j;k) \mid$  symbol $k$ occurs in cell
$(i,j)$ of the array$\}$. A {\sf partial Latin square} $P$ of
order $n$ is an $n\times n$ array with entries chosen from the set
$X = \{1,2, \ldots, n\}$, such that each element of $X$ occurs at
most once in each row and at most once in each column. Hence there
are cells in the array that may be empty, but the cells that are
filled have been filled so as to conform with the Latin property
of the array. Note that a partial Latin square of order $n$ is not
necessarily completable to a Latin square of order $n$. Let $P$ be
a partial Latin square of order $n$, then $|P|$ is said to be the
{\sf size} of the partial Latin square and the set of positions
${\cal S}_P=\{(i,j) \mid (i,j;k)\in P\}$ is said to determine the
{\sf shape} of $P$.

A partial Latin square $C$ contained in a Latin square $L$ is said
to be {\sf uniquely completable} if $L$ is the only Latin square
of order $n$ with $k$ in the cell $(i,j)$ for every $(i,j;k) \in
C$. A {\sf critical set} $C$ contained in a Latin square $L$ is a
partial Latin square that is uniquely completable, with no proper
subset of $C$ satisfying this requirement. We say a partial Latin
square $P$ {\sf forces} an entry $e=(i,j;k)$ into $P$, if $P \cup
\{ (i,j;k') \}$ is not a partial Latin square, for every $k' \neq
k$. The name ``critical set'' and the concept were invented by
statistician John Nelder, about 1977, and his ideas were first
published in a note~{\bf\cite{Nelder77}}. For a Latin square $L$,
 $\lcs{L}$ and $\scs{L}$ respectively, denote the size of the
{\sf largest critical sets} and {\sf smallest critical sets} of
$L$. Let $\lcs{n}$ be the maximum of $\lcs{L}$ over all Latin
squares $L$ of size $n$, and $\scs{n}$ be the minimum of $\scs{L}$
over all Latin squares $L$ of size $n$. Determining $\lcs{n}$ and
$\scs{n}$ are  open questions, see for example~{\bf
\cite{keedwellsurvey}}. We introduce some new bounds for
$\lcs{n}$, and for $\max \scs{L}$.

In Section 2 we show that every Latin square has a critical set of
size at most $n^2-\frac{\sqrt{\pi}}{2}n^{3/2}$, and in Section~3
we give an upper bound for the number of partial Latin squares of
order $n$ and size $k$. By using this upper bound, we prove in
Section~4 that there exist Latin squares which do not contain any
critical set of size less than $n^2-(e+o(1))n^{5/3}$. This result
improves the previously known lower bound given
in~{\bf\cite{HatamiMahmoodianLCS}}:
$$\lcs{n} \geq n^2(1-\frac{2 + \ln 2}{\ln n}) +n(1+\frac {2 {\ln
2}+\ln \left( 2 \pi \right)} {\ln n})-\frac{\ln 2}{\ln n}.$$

Note that the two bounds given in Sections 2 and 4 show that:
$$n^2-(e+o(1))n^{10/6} \le \max \scs{L} \le
n^2-\frac{\sqrt{\pi}}{2}n^{9/6}.$$

Most of our proofs are involved with calculations, and the
following well-known inequalities will be used frequently.
\begin{equation}
\label{ineq1} {a \choose b} = {a \choose a-b} \le
(\frac{ea}{b})^b,
\end{equation}
where $a$ and $b$ are natural numbers.
\begin{equation}
\label{ineq2} \sqrt{2\pi n} \left( n \over e \right) ^ n \le n!
\le \sqrt{2\pi n} \left( n \over e \right) ^ n e^{\frac{1}{12n}}.
\end{equation}

%%%%%%%%%%%%%%%%%%%%%%%%%%%%%%%%%%%%%%%%%%%%%%%%%%%%%%%%%%%%%%%%%%%%%%%%%%%%

\section{The upper bound}
In this section we use the probabilistic method to obtain an upper
bound for the size of the smallest critical set of any arbitrary
Latin square of order $n$.

\begin{theorem}
\label{upper}
 Every Latin square $L$ of order $n$ contains a
critical set of size less than $n^2-n \frac{\sqrt{n \pi}}{2}$.
\end{theorem}
\begin{proof}
Assign to each entry $e=(i,j;k)$ of $L$ a ``birth time'' $x_e$.
These $x_e$ are independent real variables, each with a uniform
distribution in $[0,1]$. Next, order the entries according to
increasing birth time, giving the ordering $e_1,e_2, \ldots,
e_{n^2}$. So we have $x_{e_i} < x_{e_{i+1}}$. Now begin from the
empty set $C$, and for every $ 1\le i \le n^2$, if the partial
Latin square $\{e_1, e_2, \ldots, e_{i-1}\}$ does not force $e_i$,
add $e_i$ to $C$. It is trivial that the constructed set $C$ is a
uniquely completable set. We want to calculate the expected size
of $C$. Consider an entry $e=(i,j;k)$ with the birth time $x_e$.
For every element $k'$ in $\{1,2,\ldots, n\} \setminus \{k\}$,
there exists a cell in the $i$-th row of $L$ and a cell in the
$j$-th column with value $k'$. Since birth times have uniform
distribution in $[0,1]$, the probability that at least one of
these two entries has birth time less than $x_e$ is $1-(1-x_e)^2$.
Thus the entry $e$ is forced by the previous entries with
probability $(1-(1-x_e)^2)^{n-1}$. So we have
$$E(|C|)=n^2(1-\int^1_0 (1-(1-x)^2)^{n-1} dx),$$
and if $\sin(\alpha)= 1-x$, then
$$E(|C|)=n^2(1-\int^{\frac{\pi}{2}}_0 \cos^{2n-1}\alpha \ d\alpha)=
n^2(1-\frac{(2n-2)(2n-4)\ldots 2}{(2n-1)(2n-3)\ldots 1}).$$

We know that $\frac{\pi}{2} \le \frac{2}{1} \frac{2}{3}
\frac{4}{3} \frac{4}{5} \ldots \frac{2n-2}{2n-1} \frac{2n}{2n-1}$
(see for example~{\bf\cite{SpiegelHandbook}}, page 188). So
$E(|C|) \le n^2-n \frac{\sqrt{n \pi}}{2}$. This implies that
there exists a uniquely completable set of size less than $n^2-n
\frac{\sqrt{n \pi}}{2}$.
\end{proof}

\section{Number of partial Latin squares}

In this section we give an upper bound for the number of partial
Latin squares of order $n$ and of size $k$. This result will be
used in Section 4. The following lemma is a corollary of
Br\'{e}gman's well-known  inequality (see for
example~{\bf\cite{vanLintWilson}}, page 83). Note that for an $m
\times n$ matrix, $A=[a_{ij}]$, $m \le n$, the permanent is
defined as
$$\per(A)=\sum_\sigma \prod_{i=1}^{m}a_{i\sigma(i)},$$
where $\sigma$ is a one to one function from $\{1,\ldots,m\}$ to
$\{1,\ldots,n\}$.

\begin{lemma}
\label{lembreg} Let $A$ be a $(0,1)$-matrix of order $m \times n$,
which has $r_i$ ones in row $i$, $1 \le i \le m.$ Then
$$\per(A) \le \frac{n!^{\frac{n-m}{n}}}{(n-m)!} \prod^{m}_{i=1}(r_i!)^{1/r_i}.
\qquad \qquad \qquad
\qquad \qquad\qquad \qquad\qquad \qquad\qquad \
{\hfill\rule{2mm}{2mm}}$$

\end{lemma}

\begin{theorem}
Let $S$ be a set of the cells in an $n \times n$ array that has
exactly $r_i$ cells in the $i$-th row and $c_j$ cells in the
$j$-th column. The number of partial Latin squares of shape $S$ is
less than or equal to
$$\left(\prod^{n}_{i=1}n!^{\frac{n-r_i}{n}}\frac{1}{(n-r_i)!}
\right)\left(\prod^{n}_{i=1}\prod^{c_i-1}_{j=0}
(n-j)!^{\frac{1}{n-j}}\right).$$
\end{theorem}
\begin{proof}
Suppose that all cells in the first $t-1$ rows of $S$ are filled,
and denote the constructed partial Latin square by $P_{t-1}$. Then
for the $t$-th row construct an $r_{t}$ by $n$ $(0,1)$-matrix
$A_t$ as in the following. For every cell $(t,y_i)$ in $S$ $(1 \le
i \le r_i)$, let the cell $(i,k)$ of $A_{t}$ be $1$ if $k$ does
not occur in the column $y_i$ of $P_{t-1}$ and $0$ otherwise. Note
that the $t$-th row of $S$ can be filled in exactly $\per(A_t)$
ways. Now if $l_{y_i}$ is the number of cells in the first $t-1$
rows of the $y_i$-th column of $S$, then by Lemma~\ref{lembreg}
$${\rm per}(A_{t}) \le \frac{n!^{\frac{n-r_{t}}{n}}}{(n-r_t)!}
\prod^{r_{t}}_{i=1}(n-l_{y_i})!^{\frac{1}{n-l_{y_i}}}.$$
 So by multiplying right sides together for $t=1,\ldots,n$, we achieve an upper bound for the number of ways that $S$ can
be filled. It is easy to see that the product is
$$\left(\prod^{n}_{i=1}n!^{\frac{n-r_i}{n}}\frac{1}{(n-r_i)!}\right)
\left(\prod^{n}_{i=1}\prod^{c_i-1}_{j=0}
(n-j)!^{\frac{1}{n-j}}\right).$$
\end{proof}\\
%\vspace*{-1cm}

\begin{theorem}
\label{number} The number of partial Latin squares of order $n$
and of size $k$ is bounded above, by:
$${n^2 \choose k}\frac{n!^{2n-\frac{k}{n}} e^{n(3+\frac{\ln(2\pi
n)^2}{4})}}{(n-\frac{k}{n})!^{2n} e^k}.$$
\end{theorem}
\begin{proof}
Let $S_P$ be a shape which has $r_i$ cells in the $i$-th row and
$c_j$ cells in the $j$-th column. Then
$r_1+r_2+\ldots+r_n=c_1+c_2+\ldots+c_n=k$. First we show that
$\prod^{n}_{i=1}n!^{\frac{n-r_i}{n}}\frac{1}{(n-r_i)!}$ achieves
its maximum value when $r_i=k/n$, for all $1\le i \le n.$ Recall
that when $x$ is any real number, $x!$ is defined as
$x!=\Gamma(x+1)$. We have
$$\prod^{n}_{i=1}n!^{\frac{n-r_i}{n}}\frac{1}{(n-r_i)!}=n!^{(n-\frac{k}{n})}
\prod^{n}_{i=1}\frac{1}{(n-r_i)!}.$$
Since $\ln \Gamma(n)$ is
convex, the expression is maximized when $r_i$ are all equal.
Hence
\begin{equation}
\label{first}
\prod^{n}_{i=1}n!^{\frac{n-r_i}{n}}\frac{1}{(n-r_i)!}\le
\frac{n!^{(n-\frac{k}{n})}}{(n-\frac{k}{n})!^n}
\end{equation}

Next consider $\prod^{n}_{i=1}\prod^{c_i-1}_{j=0}
(n-j)!^{\frac{1}{n-j}}$.  By Inequality~(\ref{ineq2}) we have
$(n-j)! \le (\frac{n-j}{e})^{n-j}\sqrt{2\pi
(n-j)}e^{\frac{1}{12(n-j)}}$, so that
$$\prod^{n}_{i=1}\prod^{c_i-1}_{j=0} (n-j)!^{\frac{1}{n-j}} \le
\prod^{n}_{i=1}\prod^{c_i-1}_{j=0} \frac{n-j}{e} (2 \pi
(n-j))^{\frac{1}{2(n-j)}}e^{\frac{1}{12(n-j)^2}}$$
$$\le \left(
\prod^n_{i=1}\prod^{c_i-1}_{j=0}\frac{n-j}{e}\right)
\left(\prod^{n}_{i=1}(2\pi
i)^{\frac{n}{2i}}e^{\frac{n}{12i^2}}\right).$$
Now knowing that $\frac{\pi^2}{6} \ge \sum_{i=1}^{\infty}
\frac{1}{i^2}$, we have
$$\ln \left(\prod^{n }_{i=1}(2\pi i)^{\frac{n}{2i}}e^{\frac{n}{12i^2}}\right)=
n\sum_{i=1}^{n}(\frac{\ln(2\pi i)}{2i}+\frac{1}{12i^2}) \le
n\left(\frac{\pi^2}{72}+\frac{\ln(2\pi)}{2}+\int_1^{n}(\frac{\ln(2\pi
x)}{2x})dx\right).$$
 The integral is equal to $\frac{\ln(2\pi n)^2}{4}$. So we have
$$\ln \left(\prod^{n }_{i=1}(2\pi i)^{\frac{n}{2i}}e^{\frac{n}{12i^2}}\right)
\le n\left(\frac{\pi^2}{72}+\frac{\ln(2\pi)}{2}+ \frac{\ln(2\pi
n)^2}{4} \right) \le n(3+\frac{\ln(2\pi n)^2}{4}).$$
For $\prod^n_{i=1}\prod^{c_i-1}_{j=0}\frac{n-j}{e},$ we have
$$ \prod^n_{i=1}\prod^{c_i-1}_{j=0}\frac{n-j}{e} = \prod_{i=1}^n\frac{n!}{(n-c_i)!e^{c_i}}=
\frac{n!^n}{e^k}\prod_{i=1}^n\frac{1}{(n-c_i)!}.$$ And again since
$\ln \Gamma(n)$ is convex, this expression achieves its maximum
value, when $c_i$ are all equal, i.e.
\begin{equation}
\label{sec} \left(\prod^{n}_{i=1}\prod^{c_i-1}_{j=0}
(n-j)!^{\frac{1}{n-j}}\right) \le \frac{n!^n e^{n(3+\frac{\ln(2\pi
n)^2}{4})} }{(n-\frac{k}{n})!^n e^k}
\end{equation}
Note that we can choose the shape of these partial Latin squares
in ${n^2 \choose k}$ ways. This fact and
Inequalities~(\ref{first}) and (\ref{sec}) lead to the result of
the theorem.
\end{proof}

\section{The lower bound}
\begin{theorem}
\label{lower} There exists a Latin square $L$ such that $\scs{L}
\ge n^2-(e+o(1))n^{5/3} $.
\end{theorem}
\begin{proof}
 As a result of van
der Waerden conjecture, we have
$$L(n) \ge \frac{(n!)^{2n}}{n^{n^2}}. $$
(see for example Theorem 17.2 in~{\bf\cite{vanLintWilson}}) where
$L(n)$ is the number of Latin squares of order $n$.

If every Latin square has a critical set of size at most $k$, then
obviously the number of critical sets of size at most $k$ is
greater than or equal to $L(n)$, and as a result, the number of
uniquely completable partial Latin squares of size $k$ is greater
than or equal to $L(n).$  By Theorem~\ref{number}, we know that
the number of partial Latin squares of size $k$ is at most
$${n^2 \choose k}\frac{n!^{2n-\frac{k}{n}} e^{n(3+\frac{\ln(2\pi
n)^2}{4})}}{(n-\frac{k}{n})!^{2n} e^k}.$$
So
$$\frac{(n!)^{2n}}{n^{n^2}} \le {n^2 \choose k}\frac{n!^{2n-\frac{k}{n}} e^{n(3+\frac{\ln(2\pi
n)^2}{4})}}{(n-\frac{k}{n})!^{2n} e^k}$$ Let \
$c=1-\frac{k}{n^2}$. \ Then by Inequality~(\ref{ineq1}), ${n^2
\choose k}={n^2 \choose n^2-k} \le ( \frac{e}{c})^{cn^2}$. Hence
$$\frac{n!^{n-cn}}{n^{n^2}} \le \frac{e^{cn^2}e^{n\ln(2\pi
n)^2}}{c^{cn^2}(cn)!^{2n}e^{n^2-cn^2}},$$ or
$$\frac{n^{n^2-cn^2}}{e^{n^2-cn^2}n^{n^2}} \le \frac{e^{3cn^2}e^{n\ln(2\pi
n)^2}}{c^{cn^2}(cn)^{2cn^2}e^{n^2-cn^2}},$$
or
\begin{equation}
\label{ineq3} c^{3c}n^c \le e^{3c}e^{\frac{\ln(2\pi n)^2}{n}}.
\end{equation}
Fix a sufficiently large $n$. If $c \ge \frac{1}{n^{1/3}}$, then
$c^{3c}n^ce^{-3c}$ increases as $c$ increases, and if
$c=\frac{e^{1+\frac{1}{\sqrt{n}}}}{n^{1/3}}$, then $c^{3c}n^c >
e^{3c}e^{\frac{\ln(2\pi n)^2}{n}}$. So Inequality~(\ref{ineq3})
implies that $c \le \frac{e+o(1)}{n^{1/3}}$, or $k \ge
n^2-(e+o(1))n^{5/3}.$
\end{proof}
\\

Cavenagh  ({\bf\cite{CavenaghThesis}}, Corollary 10.8, page 147)
proved that $\max \scs{L} < n^2-O(n^{4/3})$. Theorem~\ref{upper}
improves this result. Also Theorem~\ref{lower} shows that the a
conjecture of Cavenagh ({\bf\cite{CavenaghThesis}}, Conjecture
10.9, page 147) which states that $\max \scs{L} \le
\frac{n^2}{2}$, is not true for sufficiently large $n$.
%%%%%%%%%%%%%%%%%%%%%%%%%%%%%%%%%%%%%%%%%%55
\section*{Acknowledgements}
We thank R. Tusserkani for his valuable discussions leading
towards the result of this paper.
%%%%%%%%%%%%%%%%%%%%%%%%%%%%%%%%%%%%%%%
%\bibliographystyle{plain}
%\bibliography{sr50}

\newpage
\def\cprime{$'$}

\end{document}